\documentclass[10pt]{amsart}
\usepackage[all]{xy}
\pdfoutput=1
\usepackage{amsmath, amsthm, amssymb,slashed,stmaryrd}
\usepackage{enumitem}
\usepackage{ifpdf}
\usepackage[pdftex]{graphicx}
\usepackage{tikz}
\usetikzlibrary{matrix,arrows,calc}
\usepackage{mathtools}

\usepackage[pdftex,plainpages=false,hypertexnames=false,pdfpagelabels]{hyperref}
\usepackage[capitalise]{cleveref}

\usepackage{amssymb}
\usepackage{tikz-cd}
\usetikzlibrary{arrows,decorations.markings}
\usetikzlibrary{shapes.geometric}

\usepackage{bm}
\usepackage{accents}
\usepackage[super]{nth}
\usepackage{mathrsfs}

\usepackage{appendix}

\newtheorem{thm}{Theorem}
\newtheorem{thmalpha}{Theorem}

\newtheorem{lem}[thm]{Lemma}

\theoremstyle{definition}

\newtheorem*{setup*}{Setup}

\newtheorem{setup}[thm]{Setup}

\newtheorem{rmk}[thm]{Remark}

\newtheorem*{conjecture*}{Conjecture}
\newtheorem*{theorem*}{Theorem}
\newtheorem*{claim*}{Claim}
\newtheorem*{corollary*}{Corollary}
\newtheorem*{notation*}{Notation}

\DeclareSymbolFont{bbold}{U}{bbold}{m}{n}
\DeclareSymbolFontAlphabet{\mathbbold}{bbold}

\def\Sym{{\rm Sym}}

\def\Spec{{\rm Spec}}

\def\K3{{\rm K3}}

\def\Fil{{\rm Fil}}

\def\SO{{\rm SO}}

\def\gr{{\rm gr}}

\def\det{{\rm det}}

\newcommand{\Ql}{\mathbb{Q}_{\ell}}
\newcommand{\Qlbar}{\overline{\mathbb{Q}}_{\ell}}

\newcommand{\C}{\mathbb C}

\makeatletter
\newcommand{\colim@}[2]{%
  \vtop{\m@th\ialign{##\cr
    \hfil$#1\operator@font colim$\hfil\cr
    \noalign{\nointerlineskip\kern1.5\ex@}#2\cr
    \noalign{\nointerlineskip\kern-\ex@}\cr}}%
}
\newcommand{\colim}{%
  \mathop{\mathpalette\colim@{}}\nmlimits@
}
\makeatother

\usepackage{color}

\newcommand\nc{\newcommand}
\nc\mf\mathfrak
\nc\mc\mathcal
\nc\mb\mathbb
\nc\msf\mathsf
\nc\mscr\mathscr

\usepackage[
backend=biber,
style=alphabetic,
sorting=nyt,
maxnames=50
]{biblatex}

\usepackage[top=2cm,bottom=2cm,left=2cm,right=2cm,asymmetric]{geometry}
\begin{document}

\title[Constructing abelian varieties from rank 3 Galois representations]{Constructing abelian varieties from rank 3 Galois representations with real trace field}
\author{Raju Krishnamoorthy}

\email{krishnamoorthy@alum.mit.edu}  

\author{Yeuk Hay Joshua Lam}
\email{joshua.lam@hu-berlin.de}
\address{Humboldt Universität Berlin,
       Institut für Mathematik- Alg.Geo.,
        Rudower Chaussee 25
        Berlin, Germany}

\date{\today}

\begin{abstract}
Let $U/K$ be a smooth affine curve over a number field and let $L$ be an irreducible rank 3 $\overline{\mathbb Q}_{\ell}$-local system on $U$ with trivial determinant and infinite geometric monodromy around a cusp. Suppose further that $L$ extends to an integral model such that the Frobenius traces are contained in a fixed totally real number field. Then, after potentially shrinking $U$, there exists an abelian scheme $f\colon B_U\rightarrow U$ such that $L$ is a summand of $R^2f_*\overline{\mathbb Q}_{\ell}(1)$. 
 
The key ingredients are: (1) the totally real assumption implies $L$ admits a square root $M$; (2) the trace field of $M$ is sufficiently bounded, allowing us to use \cite{kyz} to construct an abelian scheme over $U_{\bar K}$ geometrically realizing $L$; and (3) Deligne's weight-monodromy theorem and the Rapoport-Zink spectral sequence, which allow us to pin down the arithmetizations using the total degeneration.
\end{abstract}

\maketitle 
\setcounter{tocdepth}{1}

\begin{setup}\label{setup}
    Let $X/K$ be a smooth projective geometrically connected curve over a number field, let $D\subset X$ be a reduced divisor, and set $U:=X\setminus D$. Fix a prime $\ell$ and an algebraic closure $\Qlbar$ of $\Ql$. Let $L$ be a rank $r\geq 1$ irreducible $\Qlbar$-local system on $U$ with trivial determinant.
\end{setup} It is well known that there exists an integral model $(\mf X, \mf D)$, defined over $\mc{O}:=\mc{O}_{K, S}$ for $S$ a finite set of primes of $\mc{O}_K$,  of $(X, D)$, such that  $L$ naturally extends to a $\Qlbar$-local system  $\mathscr L$ on $\mf U:=\mf{X}\setminus \mf{D}$, still with trivial determinant (see e.g. \cite[Proposition 6.1]{sasha}).

The relative Fontaine--Mazur conjecture, as reformulated by Petrov, predicts that local systems in Setup \ref{setup} occur in the cohomology of a smooth projective family after possible restricting to an open $V\subset U$ and twisting by a character of $G_K$. In this article, we make some progress towards this conjecture when the rank is 3.

\begin{thmalpha}\label{thm:main}
        Suppose $r=3$ and $L$ has infinite geometric monodromy around at least one point of $D$. Suppose further there exists a totally real number subfield $E\subset \Qlbar$ such that the field generated by Frobenius traces of $\mathscr L$ is contained in $E$. Then after potentially increasing $D$, there exists an abelian scheme $f\colon B_U\rightarrow U$ such that $L$ is a summand of $R^2f_*\Qlbar(1)$. 
\end{thmalpha}
\begin{rmk}
    Note that if we do not make any assumption on the trace field being totally real, then the conclusion of Theorem \ref{thm:main} is false. Indeed, there exist rank 3 hypergeometric local systems (with monodromy Zariski dense in $\text{SL}(3,\Qlbar)$) with bounded Frobenius trace field that do not come from abelian varieties.
\end{rmk}
An auxiliary result, which may be of indepedent interest, is the following.
\begin{lem}\label{lemma:finitedet}
        Let $U/K$ be a smooth curve over a number field with compactification $X$, let $f\colon Y\rightarrow U$ be a smooth projective morphism, and let $Q$ be a rank $r$ summand of $R^{r-1}f_*\Qlbar$. Further suppose that $Q$ is totally degenerating at some $\infty\in X\setminus U$, i.e., the local inertia around $\infty$ acts with a maximal Jordan block. Then $Q$ has   determinant isomorphic to $\Qlbar(-r(r-1)/2)$ up to a finite order character.
    \end{lem} 
    \begin{rmk}
    Our assumptions in \cref{lemma:finitedet} include  that of $Q$ being of rank $r$ and appearing in weight $r-1$. Note that there must be some relation between these quantities for the conclusion of the lemma, namely that the determinant is a power of cyclotomic (up to a finite order character), to hold. For example, for any $Q$ as in \cref{lemma:finitedet}, and $\chi$ a CM character appearing in $H^1$ of a CM abelian variety $A_0/K$, $Q\otimes \chi$ now has rank $r$ and weight $r$, and no power of the determinant is a power of cyclotomic.  
\end{rmk}
\begin{proof}[Proof of \cref{thm:main}]
We break the proof up into steps.
\begin{enumerate}[wide, labelwidth=!, labelindent=0pt]
    \item First of all, we claim that it is sufficient to prove the theorem after  replacing  $K$ by a finite extension and $U_K$ with a (geometrically connected) finite \'etale cover. Indeed, suppose we have a finite \'etale map $V\rightarrow U$ and an abelian scheme $g\colon B_V\rightarrow V$ such that $L(-1)|_V$ is a summand of $R^2g_*\Qlbar$.  Then we claim that the Weil restriction $f\colon A_U:=\mf{Res}^{V}_{U}B_V\rightarrow U$, an abelian scheme on $U$, will satisfy the conclusion of the theorem. The two relevant facts we need: $R^2g_*\Qlbar$ is semi-simple by \cite[Ch. 6, Sec. 3, Theorem I on p. 211]{rationalpoints}, and if $H\subset G$ is a finite index subgroup and $V$ a representation of $H$, then $\text{Ind}^G_H(\bigwedge^2 V)\hookrightarrow \bigwedge^2 \text{Ind}^G_H(V)$.

    \item Fix $\mf p\in \text{Spec}(\mc O)$ with residue field $\mb F_q$ and set $\mathscr L_{\mf p}:=\mathscr L|_{\mf U_{\mf p}}$. First of all, we claim that the monodromy of $\mathscr L_{\mf p}$ is contained in $\SO(3, \Qlbar)$. To see this, fix  an isomorphism $\iota\colon \Qlbar\rightarrow \C$ and set $\sigma \in \text{Aut}(\Qlbar)$ to be the transport of structure of complex conjugation. 
    
    \cite[Th\'eor\`eme VII.6]{lafforgue} implies that  every Frobenius eigenvalue of $\mathscr L_{\mf p}$ is a $q$-Weil number of weight 0 and is contained in a CM number field. Therefore the Frobenius trace field (see \cite[Definition 1.1.4]{krishnamoorthy-lam} for the definition) of $\mathscr L_{\mf p}$ is contained in a CM number field. Hence $\sigma$ acts by complex conjugation on the traces of Frobenius conjugacy classes. Let $^{\sigma}(\mathscr L_{\mf p})$ denote the $\sigma$-companion of $\mathscr L_{\mf p}$, which exists again by \cite[Th\'eor\`eme VII.6]{lafforgue}. Then $^{\sigma}(\mathscr L_{\mf p})\cong \mathscr L_{\mf p}$ by the totally real assumption. On the other hand, as $\mathscr L_{\mf p}$ is pure of weight 0, it follows that   $^{\sigma}(\mathscr L_{\mf p})\cong \mathscr L^{\vee}_{\mf p}$. Therefore there is an equivariant pairing $\mathscr L_{\mf p}\otimes \mathscr L_{\mf p}\rightarrow \Qlbar$, which is moreover non-degenerate as $\mscr{L}_{\mf p}$ is irreducible. As $\mscr{L}_{\mf p}$ has odd rank, this must be a symmetric pairing, i.e. the monodromy lies in $\SO(3, \Qlbar)$.

    \item There is an isogeny $\text{SL}(2,\Qlbar)\rightarrow \text{SO}(3,\Qlbar)$, given by the second symmetric square of the defining representation. Therefore, after possibly extending $K$ and replacing $U$ by a finite \'etale cover, we may assume $L\cong \text{Sym}^2M$, with $M$ an irreducible rank 2 $\Qlbar$-local system on $U$, with trivial determinant. Note that $M$ still has infinite monodromy around a point of $D$. Moreover, by \cite[Proposition 6.1]{sasha}, after potentially enlarging $S$, $M$ canonically extends to a $\Qlbar$ local system $\mathscr M$ on $\mf U$ with trivial determinant.

    \item Enlarge $S$ to contain 2. We claim that there exists a finite \'etale map $\mc O\rightarrow \mc O'$, a curve $\mf V/\mc O'$, and a finite \'etale morphism $\mf V\rightarrow \mf U$ such that the stable trace field (defined in \cite[Definition 1.1.4]{krishnamoorthy-lam}) of $\mathscr M|_{\mf V}$ is contained in $E$. Indeed, let $V\rightarrow U$ be a cover that trivializes the 2-torsion of the Jacobian. Let us prove that $\mathscr M|_{\mf V}$ has stable trace field contained in $E$.
    
    Fix a prime $\mf p'$ of $\mc O'$ lying over a prime $\mf p$ of $\mc O$, let the trace field of $\mathscr M_{\mf p}$ be $F\supset E$, and consider the embeddings $\mc S:=\{\sigma\colon F\hookrightarrow \Qlbar\}$ that fix $E$ elementwise. Then, for each $\sigma\in \mc{S}$, the companion   $^{\sigma}(\mathscr M_{\mf p})$ has  second symmetric square isomorphic to $L$, as the operation of companions commutes with tensorial constructions \cite[Proof of Lemma 2.5]{krishnapal2}. Hence for each $\sigma\in \mc S$, there exists a 2-torsion rank 1 local system $R$ on $\mf U_{\mf p}$ such that  $^{\sigma}(\mathscr M_{\mf p})\cong \mathscr M_{\mf p}\otimes R$. Therefore we have $^{\sigma }(\mathscr M_{\mf p})|_{\mf V_{\mf p'}}\cong \mathscr M_{\mf p}|_{\mf V_{\mf p'}}\otimes \psi$, with $\psi$ being a 2-torsion rank 1 local system of the base $\Spec(\mc O'/\mf p')$. 
    
    Therefore, all of the local systems $^{\sigma} (\mathscr M_{\mf p})|_{\mf V_{\mf p'}}$ become isomorphic after a quadratic extension of the base $\mc O'/\mf p'$, which implies that the stable trace field of $\mathscr M|_{\mf V}$ is contained in $E$. The upshot is that,  replacing  $U$ by $V$, we may assume that the stable trace field of $M$ is bounded.
    \item The boundedness of the stable trace field of $M$ implies, using the argument of \cite{kyz} (using the moduli space $\mc H$, not the more refined moduli spaces $\mc H_{i}$ or $\mc H_{\infty}$ of \emph{loc. cit.}) that there exists a principally polarized abelian scheme $\bar{f}\colon A_{U_{\bar K}}\rightarrow U_{\bar K}$ such that $R^1\bar{f}_*\Qlbar$ has $M|_{U_{\bar K}}$ as a summand.\footnote{This argument is also indicated in Remark 1.9 of \emph{loc. cit.} Here is a sketch. After inverting finitely many primes, $\mathcal H/\mc O$ of Definition 3.1 of \emph{loc. cit.} is finite flat by Section 4. Moreover, $\mathcal H$ and has mod $p$ points for infinitely many $p$ by Section 2; hence has a characteristic 0 point. Then the prime-to-$p$ specialization isomorphism of $\pi_1$ implies the desired result. See also Footnote 5 of \emph{loc. cit.}} As the abelian scheme $A_{U_{\bar{K}}}\rightarrow U_{\bar K}$ is rigid (see Section 4 of \emph{loc. cit.}), it follows from \cite[Theorem 2]{faltings} that the summand $M|_{U_{\bar K}}$ is cut out by an element of $End(A_{U_{\bar K}})\otimes \Qlbar$. Therefore, after replacing $K$ by a finite extension (using that the endomorphism ring of an abelian scheme is a finitely free $\mb Z$-module), it follows that there exists an abelian scheme $f: A\rightarrow U$ such that the cohomology has a rank 2 factor $M'$ with the following property: $M_{U_{\bar K}}\cong M'_{U_{\bar K}}$.   It now suffices to show that there exists a finite  \'etale cover $V\rightarrow U$ such that $M|_V\simeq M'|_V(1/2)$ for some choice of  $\Qlbar(1/2)$: indeed, this implies that $L|_V\simeq (\Sym^2M)|_V\simeq \Sym^2(M'|_V)(1)$, and by construction the latter  appears in $R^2f_*\Qlbar(1)$.  

    \item By \cref{lemma:finitedet}, $M'$ has determinant equal to $\Qlbar(-1)$ up to a finite order character, and we conclude that 
    \[M\simeq M'(1/2)\otimes \chi\]
    for any choice of $\Qlbar(1/2)$, and $\chi$ some  finite order rank one local system on $U$. Passing to the finite cover  $V\rightarrow U$ which trivializes $\chi$ gives $M|_V\simeq M'|_V$ as desired.
\end{enumerate}
\end{proof}

    \begin{proof}[Proof of \cref{lemma:finitedet}]
        
        Since the conclusion of our lemma is insensitive to passing to finite covers of $U$, we may assume that $Q$  has maximal unipotent degeneration around some cusp and moreover the map $f\colon Y\rightarrow U$ has strict semistable reduction.
        
       By Cebotarev density, it suffices to prove the claim  after reduction mod $\mf p$ for almost all primes \footnote{using also the fact that, for an $\ell$-adic local field $F$ and a number field $K$,  an almost everywhere unramified character $\lambda: G_K\rightarrow \mc{O}_F^{\times}$ is of finite order iff the Frobenius elements at almost all primes have finite (though a priori unbounded) orders}.  Let $\infty$ be a cusp of $U_{\mf{p}}$, around which $Q$ is totally degenerating, with uniformizing parameter $z_{\infty}$. It suffices to address the formal local picture around $k((z_{\infty}))$, where $k$ is the field of definition of $\infty$. More precisely, if $V$ is the $G_{k((z_{\infty}))}$-representation corresponding to $Q|_{\Spec(k((z_{\infty})))}$, we will show that $\det V$ is  isomorphic to $\Qlbar(-r(r-1)/2)$, up to a finite order character.
        
        The maximal unipotent degeneration assumption implies the monodromy filtration $\Fil^N_i$ (see e.g. \cite[Definition/Proposition 9.2]{scholze} for a definition) has associated graded satisfying 
        \[
        \dim_{\Qlbar} \gr^N_{-(r-1)+2k} =1, \dim_{\Qlbar} \gr^N_{-(r-1)+2k+1}=0
        \]
        for $k=0, \cdots , r-1$; furthermore, the monodromy operator
        \begin{equation}\label{eqn:monooperator}
        N: \gr^N_{-(r-1)+2k+2}\rightarrow \gr^N_{-(r-1)+2k}
        \end{equation}
        is a $G_{k}$-equivariant  isomorphism. 
        
        Moreover $Q^{I} = \Fil^N_{-(r-1)}$ is  one dimensional, where $I\subset G_{k((z_{\infty}))}$ denotes the inertia subgroup. Let $\Phi\in G_k$ denote a geometric Frobenius element. The $\ell$-adic Steenbrink spectral sequence of Rapoport-Zink (\cite[Satz 2.10]{rapoport-zink}) implies that  $\Phi$ acts on  $Q^I$  via an algebraic integer; on the other hand, the weight-monodromy conjecture, in this setting a theorem of Deligne \cite[Corollaire 1.8.5]{weilii}, implies that it is pure of weight 0, and we deduce that  $\Phi$ acts on $Q^I$ by a root of unity $\zeta$. By Equation \eqref{eqn:monooperator} and the commutation relation $N\Phi=q\Phi N$, the $\Phi$-action on $\gr^N_{-(r-1)+2k}$ is given by $q^{k}\zeta$, and hence the $\Phi$-action on $\det V$ is given by $\zeta^{r}q^{1+\cdots +(r-1)}=\zeta^r q^{\frac{(r-1)r}{2}}$, as desired.
    \end{proof}

\begin{rmk}
    The hypothesis of Theorem \ref{thm:main} may be weakend to the following: the Frobenius trace field of $L$ is bounded and there exists a prime $\mf p$ of $\mc O$ such that the Frobenius trace field of $\mathscr L_{\mf p}$ is totally real. Indeed, the hypothesis will imply that the geometric monodromy of $L$ lands in $\text{SO}(3, \Qlbar)\subset \text{SL}(3,\Qlbar)$. Then we claim that for any other prime $\mf p'$, the arithmetic monodromy of $\mathscr L_{\mf p'}$ must also be orthogonal; indeed, this follows from the fact that the geometric monodromy group is a normal subgroup of the arithmetic monodromy group and that the group $\text{SL}(3,\Qlbar)$ is almost simple.
\end{rmk}

\begin{rmk}
    Step (5) of the proof of \cref{thm:main}, and in particular \cref{lemma:finitedet} obviates the use of the more delicate moduli spaces $\mathcal H_k$ in \cite[Proof of Theorem 1.4]{kyz}, c.f. Definition 3.5 and Lemma 3.6 of \emph{loc. cit.}
\end{rmk}
\subsection*{Acknowledgments} During the course of this work, Lam was supported by a Dirichlet Fellowship and Krishnamoorthy was supported by the European Research Council (ERC)
under the European Union’s Horizon 2020 research and innovation program, grant agreement
no. 101020009, project TameHodge.

\printbibliography[]

\end{document}